\documentclass[12pt]{amsart}
\usepackage{amsmath,amscd,amssymb,amsfonts}

\newcommand{\bC}{{\mathbf C}}

\newcommand{\bP}{{\mathbf P}}

\newcommand{\bQ}{{\mathbf Q}}
\newcommand{\bZ}{{\mathbf Z}}

\newcommand{\cA}{{\mathcal A}}

\newcommand{\cO}{{\mathcal O}}

\newcommand{\codim}{\hbox{\rm codim}\,}

\newcommand{\ra}{\rightarrow}

\def\eps{\epsilon}

\theoremstyle{plain}
\newtheorem{thm}{Theorem}[section]

\newtheorem{conj}[thm]{Conjecture}
\newtheorem{lem}[thm]{Lemma}
\newtheorem{prop}[thm]{Proposition}

\theoremstyle{definition}

\newtheorem{rem}[thm]{Remark}
\newtheorem{example}[thm]{Example}
\newtheorem{subs}[thm]{}

\title[The Monodromy Conjecture for Hyperplane Arrangements]{The Monodromy Conjecture for Hyperplane Arrangements}
\author{Nero Budur}
\address{Department of Mathematics,
University of Notre Dame, 255 Hurley Hall, IN 46556, USA} \email{nbudur@nd.edu}
\author{Mircea Musta\c{t}\u{a}}
\address{2844 East Hall,
Department of Mathematics,
University of Michigan,
530 Church Street,
Ann Arbor, MI 48109-1043} \email{mmustata@umich.edu}
\author{Zach Teitler}\address{Department of Mathematics,
Mailstop 3368,
Texas A\&M University,
College Station, TX 77843-3368}\email{zteitler@tamu.edu}

\thanks {The first author was partially supported by the NSF grant DMS-0700360,
and the second author was partially supported by NSF grant DMS-0758454 and by a Packard
Fellowship}

\begin{document}

\begin{abstract}
The Monodromy Conjecture asserts that if $c$ is a pole of the local topological zeta function of a hypersurface, then ${\rm exp}(2\pi i c)$ is an eigenvalue of the monodromy on the cohomology of the Milnor fiber. A stronger version of the conjecture asserts that every such $c$ is a root of the Bernstein-Sato polynomial of the hypersurface. In this note we prove the weak version of the conjecture for
hyperplane arrangements. Furthermore, we reduce the strong version to the following conjecture:
 $-n/d$ is always a root of the Bernstein-Sato polynomial of an indecomposable essential central hyperplane arrangement of $d$ hyperplanes in $\bC^n$.
\end{abstract}

\maketitle

\section{Introduction}

Let $f\in\bC[x_1,\ldots ,x_n]$ be a non-constant polynomial, with $f(0)=0$.
Consider a log resolution $\mu: Y\ra \bC^n$
of $f^{-1}(0)$, and let $E_i$ for $i\in S$ be the
irreducible components of $E=(f\circ\mu)^{-1}(0)$. Denote by $a_i$
the order of vanishing of $f\circ\mu$ along $E_i$, and by $k_i$
the order of vanishing of ${\rm det}\left({\rm Jac}(\mu)\right)$
along $E_i$. The {\it local topological zeta function of} $f$ is
\begin{equation}\label{eqZeta}
Z_{top, f}(s):=\sum_{I\subseteq S}\chi(E_I^{\circ}\cap \mu^{-1}(0))\cdot\prod_{i\in I}\frac{1}{a_is+k_i+1},
\end{equation}
where $E_I^{\circ}=(\cap_{i\in I}E_i)-\cup_{i\not\in I}E_i$, and $\chi$ is the Euler-Poincar\'{e} characteristic. J.~Denef and F.~Loeser showed in \cite{DL} that $Z_{top, f}(s)$ is independent of the choice of the log resolution $\mu$. The poles of $Z_{top,f}(s)$ are among the negative rational numbers
$\{-(k_i+1)/a_i\ |\ i\in S \}$.

\medskip

The Topological Monodromy Conjecture of \cite{DL} is a
variant of the original Monodromy Conjecture of J.-i. Igusa.
To state it, we fix some notation. For $x\in f^{-1}(0)$, let
$M_{f,x}$ denote the {\it Milnor fiber} of $f$ at $x$, which is
defined as the intersection of $f^{-1}(t)$ with a small ball of
radius $\eps$ around $x$ ($0<t\ll \eps\ll 1$). As a
$C^\infty$-manifold, $M_{f,x}$ does not depend on $t$  and $\eps$.

The \emph{Bernstein-Sato polynomial} (also called \emph{$b$-function}) of $f$ is the non-zero monic polynomial
$b_f(s)\in\bC[s]$ of minimal degree among polynomials $b(s)$
satisfying
$$
b(s)f^s=P\bullet f^{s+1}
$$
for some algebraic differential operator $P\in\bC[x_1,\ldots, x_n,\frac{\partial}{\partial {x_1}},\ldots ,\frac{\partial}{\partial {x_n}},s]$.

\begin{conj}\label{conjMonConj} ${\rm (}$\cite{DL}, Topological Monodromy Conjecture${\rm )}$ If $c$ is a pole of $Z_{top,f}(s)$, then
\begin{enumerate}
\item[(a)] ${\rm exp}(2\pi ic)$ is an eigenvalue of the monodromy action on $H^i(M_{f,x},\bC)$, for some $i$ and some $x\in f^{-1}(0)$.
\item[(b)] $b_f(c)=0$. \footnote{A stronger conjecture is believed to be true, where $x$ is close to $0$, and the local $b$-function at $0$ is used. However, in our case this makes no difference as we deal with homogeneous polynomials.}
\end{enumerate}
\end{conj}
It is well-known, by results of B.~Malgrange \cite{Mal} and M.~Kashiwara \cite{Kashiwara},
that all the roots of the $b$-function
are negative rational numbers. Also, the set of numbers ${\rm exp}(2\pi ic)$, when $c$ varies over the roots of the
$b$-function of $f$, coincides with the set of eigenvalues of
the monodromy action on $H^i(M_{f,x},\bC)$ for some $i$, and some
$x\in f^{-1}(0)$. Hence part (b) of the above conjecture implies
part (a). The main result of this note is to prove part (a) in the
case of hyperplane arrangements, and to reduce part (b) in this
case to the following

\begin{conj}\label{conjBernstein}
If $g\in {\mathbf C}[x_1,\ldots,x_m]$ is a polynomial of degree $d$, such that
$\left(g^{-1}(0)\right)_{red}$ is an indecomposable, essential, central hyperplane arrangement,
then $-\frac{m}{d}$ is a root of the Bernstein-Sato polynomial $b_g(s)$.
\end{conj}

The following is our main result.

\begin{thm}\label{thmMainThm} Let $f\in {\mathbf C}[x_1,\ldots,x_n]$
define a hyperplane arrangement in $\bC^n$.
\begin{enumerate}
\item[(a)] Conjecture~\ref{conjMonConj}-(a) holds for $f$.
\item[(b)] If Conjecture~\ref{conjBernstein} holds (for all $g$), then
Conjecture~\ref{conjMonConj}-(b) holds for $f$.
\end{enumerate}
\end{thm}
In fact, in part (b) the condition has to be checked only for the
set of arrangements $f_W$ induced on the so-called dense edges $W$
of the given hyperplane arrangements, see Remark
\ref{rmk_main_theorem}. Theorem \ref{thmMainThm} (a) implies that
if $c$ is a pole of $Z_{top, f}(s)$ then $b_f(c)=0$ or
$b_f(c-1)=0$, since the roots of $b_f(s)$ are in $(-2,0)$ in this
case, by \cite{Sa}.

The original Monodromy Conjecture of J.-i.~Igusa involves $p$-adic
local zeta functions. Theorem \ref{thmMainThm} has an analogous
version for $p$-adic local zeta functions, see Theorem
\ref{thmMainThmPadic} below.

In the last section we discuss several examples, treating in
particular generic arrangements, as well as a large class of
arrangements.

Essentially all ingredients of the proof of Theorem
\ref{thmMainThm} have been worked out in the literature. We are
merely performing the task of assembling them together and
recording a result that we could not find in the literature. We
would like to thank W.~Veys and U.~Walther whose interest prompted
us to write down this note, and to M.~Saito for useful comments
and the remarks about moderate arrangements in the last section.

\section{Proof of Theorem \ref{thmMainThm}}

Let $V_j\subset\bC^n$ be the irreducible components of the reduced
 hyperplane arrangement ${\mathcal A}:=(f^{-1}(0))_{red}$.
We assume the arrangement to be \emph{central}, that is, each
$V_j$ is a linear subspace of codimension one of $\bC^n$. In this
case $f$ is homogeneous. We deal at the end of this section with
the non-central case.

 The following definitions only depend on $f_{red}$.
 The central hyperplane arrangement ${\mathcal A}$ is {\it
indecomposable} if  there is no linear change of coordinates
on $\bC^n$ such that $f$ can be written as the product of two
nonconstant polynomials in disjoint sets of variables.
An {\it edge} of ${\mathcal A}$ is
an intersection of hyperplanes $V_i$.
The arrangement ${\mathcal A}$ is \emph{essential} if $\{0\}$ is an edge of ${\mathcal A}$.
An edge $W\subset\bC^n$ is
called {\it dense} if the hyperplane arrangement ${\mathcal A}_{W}$ given by
the image of $\cup_{V_j\supset W}V_j$ in $\bC^n/W$ is
indecomposable. For example, $V_j$ is a dense edge for every $j$. Let $S$ be the set of dense edges of ${\mathcal A}$. We take
$\mu:Y\ra\bC^n$ to be obtained by successive blowups,
by taking for $d=0,1,\ldots,n-2$ the blowup along
the (proper transform of) the union of the dense
edges of dimension $d$.

\begin{prop}\label{propDP} ${\rm (}$\cite{STV} -Theorem 3.1${\rm )}$
The morphism $\mu:Y\ra \bC^n$ is a log resolution of $f^{-1}(0)$.
\end{prop}

This log resolution is the same as the wonderful model of
\cite{DP} when the building set is taken to be the minimal one,
which is exactly the set of dense edges. However, we do not need
this additional description. Let $E_W\subset Y$ be the proper transform of
the exceptional divisor corresponding to $W\in S$. We have:
$k_W=\codim W -1$ and $a_W=\sum_{W\subseteq V_j} m_j$ for all $W\in
S$, where $a_W, k_W$ are as in (\ref{eqZeta}), and $m_j$ is the
order of vanishing of $f$ along $V_j$ (hence $m_j=1$ if $f$ is reduced).

\begin{lem}\label{lemIndecChi} ${\rm (}$\cite{STV} -Proposition 2.6${\rm )}$ Let
$g\in\bC[x_1,\ldots , x_m]$ be such that ${\mathcal
B}=(g^{-1}(0))_{red}$ is a central hyperplane arrangement in
$\bC^m$. If $U=\bP^{m-1}\smallsetminus \bP (\mathcal{B})$, where
$\bP(\mathcal{B})$ is the projectivization of $\mathcal{B}$, then
${\mathcal B}$ is indecomposable if and only if $\chi (U)\ne 0$.
\end{lem}

\begin{lem}\label{lemIndecMondromy} If $g\in\bC[x_1,\ldots , x_m]$ with $\deg g=d$ gives an indecomposable central hyperplane arrangement
${\mathcal B}=(g^{-1}(0))_{red}\subset\bC^m$, then $\{{\rm
exp}(2\pi i k/d)\mid k=1, \ldots, d \}$ is the union of the sets
of eigenvalues of the monodromy action on $H^i(M_{g,0},\bC)$ with
$i=0,\ldots, m-1$.
\end{lem}
\begin{proof} It is known that for every
homogeneous polynomial $g$, the Milnor fiber $M_{g,0}$ is
diffeomorphic with $g^{-1}(1)$, and gives a finite $d$ to $1$
unramified cover of $U:=\bP^{m-1}\smallsetminus\bP({\mathcal B})$
(see \cite{Di} p.72 (1.13)). It is also well-known that each
eigenvalue of the action of the monodromy on the cohomology of the
Milnor fiber is of the form ${\rm exp}(2\pi i k/d)$ (see
\cite{CS}). Conversely, recall that the monodromy zeta function of
$g$ is by definition
$$
\prod_{0\le j\le m-1} {\rm det} ({\rm Id}-t\cdot M_j)^{(-1)^j},
$$
where $M_j$ is the monodromy action on $H^j(M_{g,0},\bC)$. By
Example 6.1.10 in \cite{D}, the monodromy zeta function of $g$ is
$(1-t^d)^{\chi(M_{g,0})/d}$.  Let $\chi_c(-)$ denote the Euler
characteristic with compact supports. Since
$\chi(M_{g,0})=\chi_c(M_{g,0})$ and $\chi(U)=\chi_c(U)$ (see for
example \cite{F},  p.141, note 13), it follows from the
multiplicativity property of $\chi_c(-)$ with respect to finite
unramified coverings that $\chi(M_{g,0})/d=\chi(U)$. It follows
from Lemma \ref{lemIndecChi} that for every $k$, ${\rm exp}(2\pi
ik/d)$ is an eigenvalue of $M_j$ on $H^j(M_{g,0},\bC^n)$, for some
$j$.
\end{proof}

We can now prove the result stated in the Introduction for the
case of central arrangements.

\begin{proof}[Proof of part (a) of Theorem \ref{thmMainThm}]
By the definition (\ref{eqZeta}) and Proposition \ref{propDP},  the
poles of $Z_{top,f}(s)$ are included in the set $\{-(\codim W)/a_W\ |\ W\in S\}$.
If $W\in S$, then ${\mathcal A}_W$ is indecomposable and
is the reduced zero locus of a product $f_W$ of linear forms on
$\bC^n/W$, with $\deg(f_W)=a_W$. By Lemma \ref{lemIndecMondromy},
$\exp({-2\pi i(\codim W )/a_W})$ is an eigenvalue of the monodromy
on the Milnor cohomology of $f_W$ at the origin in $\bC^n/W$. Now,
take a point $p\in W-\cup_{V_j\not\supset W}V_j$.
After choosing a splitting of $W\subset \bC^n$, we have locally around $p$,
${\mathcal A}={\mathcal A}_W\times W\subset \bC^n=\bC^n/W\times W$ and
$f=f_W\cdot u$, where $u$ is a (locally) invertible function.
Hence the Milnor fiber of $f_W$ at the origin is a deformation
retract of $M_{f,p}$. Therefore ${\rm exp}({-2\pi i(\codim W )/a_W})$ is an
eigenvalue of the monodromy on the cohomology of the Milnor fiber of $f$ at
$p$.
\end{proof}




\medskip

\begin{proof}[Proof of part (b) of Theorem \ref{thmMainThm}]
We assume that for all $m,d>0$, and for all polynomials $g$ of
degree $d$ in $m$ variables defining an indecomposable, essential, central
hyperplane arrangement, we have $b_g(-m/d)=0$. By Proposition
\ref{propDP}, the candidate poles of $Z_{top, f}(s)$ are $-(\codim
W)/a_W$, with $W\in S$. Let $f_W$ be induced by $f$ as before,
such that ${\mathcal A}_W=(f_W^{-1}(0))_{red}$ in $\bC^n/W$. Then
$(\codim W)/a_W=\dim (\bC^n/W)/\deg f_W$, with ${\mathcal A}_W$
indecomposable (and automatically essential and central) in $\bC^n/W$. By assumption, we have
$b_{f_W}(-(\codim W)/a_W)=0$. On the other hand, $b_{f_W}(s)$
equals the local $b$-function $b_{f,p}(s)$ for $p$ as in the proof
of part (a). But $b_f(s)$ is the least common multiple of the
local $b$-functions. Hence every candidate pole as above is a root
of $b_f(s)$.
\end{proof}

\begin{rem}\label{rmk_main_theorem}
It follows from the above proof that in order for Conjecture~\ref{conjMonConj}-(b)
to hold for $f$, it is enough to have Conjecture~\ref{conjBernstein} satisfied by all $f_W$,
with $W\in S$.
\end{rem}

\begin{rem}\label{rmkNonCentralCase}
(Non-central arrangements.) Let $f\in {\mathbf C}[x_1,\ldots,x_n]$
be such that $(f^{-1}(0))_{red}=\cup_jV_j$ is a not necessarily
central hyperplane arrangement. Let $f_0$ be the polynomial
obtained from $f$ by only considering the factors that vanish at
$0$. In this case $Z_{top, f}(s)=Z_{top,f_0}(s)$. We showed that
the exponentials of the poles of $Z_{top,f_0}(s)$ are included in
the union over indecomposables $W$ for $f_0$ of the eigenvalues of
monodromy on $M_{f_0,x}$, where $x$ is any point of $W-\cup_{0\in
V_j\not\supset W}V_j$. Also, assuming Conjecture
\ref{conjBernstein}, the poles of $Z_{top,f_0}(s)$ are included in
the union over indecomposables $W$ for $f_0$ of the roots of
$b_{f_0,x}$. We can choose the point $x$ such that $x\not\in
\cup_{0\not\in V_j}V_j$. In this case $M_{f_0,x}=M_{f,x}$ and
$b_{f_0,x}=b_{f,x}$. This shows that Theorem \ref{thmMainThm} is
valid for the case of non-central arrangements as well.
\end{rem}

\begin{subs} {\bf The $p$-adic case.} We next consider the $p$-adic version of the Monodromy
Conjecture. For a polynomial $f\in\bQ[x_1,\ldots ,x_n]$ and a
prime number $p$, Igusa's {\it $p$-adic local zeta function} of
$f$ is defined as $Z_f^{p}(s)=\int_{(\bZ_p)^n}|f|_p^s\;d\mu$,
where $d\mu$ is the Haar measure on $(\bQ_p)^n$, and $|.|_p$ is
the $p$-adic norm. $Z_f^p(s)$ is a meromorphic function for
$s\in\bC$ whose poles encode asymptotic behavior of the numbers
$N_m$ of solutions of $f$ modulo $p^m$ when $f$ is defined over
$\bZ$, see \cite{Den}. Let $Re(c)$ denote the real part of a
complex number $c$.
\end{subs}

\begin{conj}\label{conjMonConjPadic} ${\rm (}$\cite{Ig}, Monodromy Conjecture${\rm
)}$ Let $f\in\bQ[x_1,\ldots ,x_n]$. For almost all prime numbers
$p$, if $c$ is a pole of $Z_{f}^{p}(s)$, then
\begin{enumerate}
\item[(a)] ${\rm exp}(2\pi i \cdot Re(c))$ is an eigenvalue of the monodromy action on $H^i(M_{f,x},\bC)$,
for some $i$ and some $x\in\bC^n$ with $f(x)=0$.
\item[(b)] $b_f(Re(c))=0$.
\end{enumerate}
\end{conj}
\noindent It is speculated that the conjectures might be true for
all prime numbers $p$, see \cite{Den} -2.3.

The analog of Theorem \ref{thmMainThm} involves {\it
$\bQ$-hyperplane arrangements}. These are hyperplane arrangements
in $\bC^n$ such that each hyperplane, i.e. irreducible component,
is defined over $\bQ$. Note that this implies that every edge is
defined over $\bQ$. Let $f\in\bQ[x_1,\ldots ,x_n]$. We say that
$f$ defines a $\bQ$-hyperplane arrangement if $f$, viewed as a
polynomial with complex coefficients, does. In other words, $f$
splits completely into linear factors defined over $\bQ$.

\begin{thm}\label{thmMainThmPadic} Let $f\in {\mathbf Q}[x_1,\ldots,x_n]$
define a $\bQ-$hyperplane arrangement.
\begin{enumerate}
\item[(a)] Conjecture~\ref{conjMonConjPadic}-(a) holds for $f$ and all prime numbers $p$.
\item[(b)] If Conjecture~\ref{conjBernstein} holds for all $f_W$ where $W$ is a dense edge of $f$, then
Conjecture~\ref{conjMonConjPadic}-(b) holds for $f$ and all prime
numbers $p$.
\end{enumerate}
\end{thm}
\begin{proof} Let $\mu$ be the log resolution of $f$ from Proposition
\ref{propDP}. Consider the hyperplane arrangement $f^{(p)}$
defined by $f$ in $(\bQ_p)^n$. Any dense edge of $f^{(p)}$ is also
defined over $\bQ_p$. By taking successive blowing ups of
$(\bQ_p)^n$ along (proper transforms) of dense edges as in
Proposition \ref{propDP}, the map $\mu^{(p)}$ thus obtained is a
log resolution of $f^{(p)}$ over $\bQ_p$. There is a bijection
between the dense edges of $f$ and those of $f^{(p)}$, and the
orders of vanishing of $f\circ\mu$ and $f^{(p)}\circ\mu^{(p)}$
along the irreducible divisors corresponding to dense edges are
the same. Similarly for the orders of vanishing of ${\rm
det}\left({\rm Jac}(\mu)\right)$ and ${\rm det}\left({\rm
Jac}(\mu^{(p)})\right)$. Then, by \cite{Ig} -Theorem 8.2.1 (see
also bottom of p.123 of \cite{Ig}), the real parts of the poles of
$Z_f^{p}(s)$ are among the set $\{-(\codim W)/a_W\ |\ W\text{
dense edge of }f\}$. From now on, the proof is the same as for
Theorem \ref{thmMainThm}.
\end{proof}

\begin{rem} Conjecture \ref{conjMonConjPadic} was
originally stated with $\bQ$ replaced by any number field
$F\subset\bC$, $\bQ_p$ replaced by a $p$-adic completion $K$ of
$F$, and $Z^p_f(s)$ replaced by an integral over $(\cO_K)^n$, see
\cite{Den} -2.3. We get the corresponding  more general version of
Theorem \ref{thmMainThmPadic}, and the proof, by considering
$F$-hyperplane arrangements and all $p$-adic completions $K$.
\end{rem}

\section{Examples}

In this section we discuss some examples in which we can prove Conjecture~\ref{conjBernstein}.
We start with the case of generic hyperplane arrangements. Recall that a central hyperplane arrangement ${\mathcal A}$ is \emph{generic} if for every nonzero edge $W$, there are precisely
$\codim(W)$ hyperplanes in ${\mathcal A}$ containing $W$. Such an arrangement is essential
if and only if the number of hyperplanes is at least the dimension of the ambient space.

\begin{example}\label{generic}
Let $g\in\bC[x_1,\ldots,x_n]$ be a reduced polynomial, defining a
generic central, essential arrangement of $d$ hyperplanes in
$\bC^n$. It was shown in \cite{Wal} that the roots of $b_g$ are
those $-\frac{i+n}{d}$, with $0\leq i\leq 2d-n-2$. In particular,
we see that for $i=0$ we get that $-\frac{n}{d}$ is a root of
$b_g$, confirming Conjecture~\ref{conjBernstein} in this case.
(Note that a generic arrangement is decomposable if and only if
$d=n$, when in suitable coordinates the arrangement consists of
the coordinate hyperplanes.)
\end{example}

\medskip

Note that every central arrangement in $\bC^2$ is generic, hence
Example~\ref{generic} applies.

\begin{rem}\label{rem_dim3}
Let $g\in k[x,y,z]$ be a reduced polynomial, defining a central, essential arrangement of
$d$ hyperplanes in $\bC^3$.  In this case Conjecture~\ref{conjMonConj}-(b) holds for $g$ if
Conjecture~\ref{conjBernstein} holds for $g$, that is,
if either $\cA$ is decomposable, or if
$b_g(-3/d)=0$. Indeed, this follows
from Remark~\ref{rmk_main_theorem}, together with the fact
that if $W\neq\{0\}$, then $g_W$ satisfies Conjecture~\ref{conjBernstein} by Example~\ref{generic}.
\end{rem}
A recent computation, that will appear elsewhere, of M.~Saito, S.~Yuzvinsky,
and the first author, shows that Conjecture~\ref{conjBernstein} is true for $m=3$ and $g$ reduced,
and hence Conjecture~\ref{conjMonConj}-(b) is also true.

\bigskip
It is known that every jumping number in $(0,1]$ is the negative of a root of the
Bernstein-Sato polynomial, see \cite{ELSV}, Theorem B (but not conversely).
So to show Conjecture~\ref{conjBernstein} it is enough to show $m/d$ is a jumping number.
We now consider some cases in which it can be shown $m/d$ is a jumping number,
and an example in which it is not.

\medskip
We first consider reduced polynomials
defining central, essential hyperplane arrangements in $\bC^3$.
Let $g\in\bC[x,y,z]$ be such a polynomial defining the arrangement
${\mathcal A}$. We denote by $\bP(\cA)$ the corresponding
arrangement of lines in $\bP^2$. For every $m\geq 2$, let $\nu_m$
be the number of points with multiplicity $m$ in $\bP(\cA)$. Note
that $\cA$ is generic if and only if there are no points in
$\bP(\cA)$ of multiplicity $>2$, and in this case
Conjecture~\ref{conjBernstein} is satisfied by
Example~\ref{generic}. We also point out that since the
arrangement is in $\bC^3$, it is decomposable if and only if all
but one of the lines in $\bP(\cA)$ pass through a point.

\begin{example}\label{exdim3}
Suppose that $g$ is as in the previous remark.
Corollary~1 in \cite{BS} implies that
$3/d$ is a jumping number for the multiplier ideals of $g$
if and only if
\begin{equation}\label{eqSpectrum}
\sum_{2d/3<m\le d}\nu_m\neq 1
\end{equation}
(here we assume for simplicity that $d\ge 5$).
 It follows that
whenever (\ref{eqSpectrum}) holds, Conjecture~\ref{conjBernstein}
holds for $g$. In particular, using Remark~\ref{rem_dim3} and
Example~\ref{generic}, one can check for example that Conjecture
~\ref{conjMonConj}-(b) holds for all (reduced) central hyperplane
arrangements $\cA$ in $\bC^3$ such that the multiplicities at
points of $\bP(\cA)$ are $\leq 3$. We will state a stronger
result, for moderate arrangements, in Theorem \ref{corModerate}.
\end{example}

\begin{rem}\label{remCounterex}
Even in the case of hyperplane arrangements defined by a reduced
polynomial $g$, in general it is not the case that one can prove
Conjecture~\ref{conjBernstein} by showing that $m/d$ is a
jumping number for the multiplier ideals of $g$. For example, if
$g=xy(x-y)(x+y)(x+z)$, then $3/5$ is not a jumping number (one can
use, for example, Corollary~1 in \cite{BS}). However, this
arrangement is decomposable, so Conjecture \ref{conjMonConj}-(b)
holds by Remark \ref{rem_dim3}. In order to get an indecomposable
example, we consider an arrangement with $d=10$, $\nu_7=1$,
$\nu_3=3$,  and with $\nu_{m'}=0$ for $7\ne m'> 3$. In this case
we see by \emph{loc. cit.} that $3/10$ is not a jumping number.
For example, we may take
$$g=xy(x+y)(x-y)(x+2y)(x+3y)(x+4y)(2y+z)(x+2y+z)z.$$
On the other hand, it can be checked that
$-3/10$ is indeed a root of $b_g(s)$, as predicted by Conjecture \ref{conjBernstein}. For this one uses Theorem 4.2 (e) in \cite{Sa}, with $k=3$, and $I=\{1,2\}$ as a subset of the set $\{1,\ldots , 10\}$ that remembers the factors of $g$ in the order written.
\end{rem}

\medskip

The following remarks were pointed out to us by M.~Saito. Let
$\cA$ be a hyperplane arrangement. For an edge $W$, let
$a_W$ be, as before, the sum $\sum_{W \subseteq V_j} m_j$
where $m_j$ is the multiplicity of the hyperplane $V_j$ in $\cA$.
We say that $\cA$ is a hyperplane arrangement {\it of
moderate type} if the following condition is satisfied:
\begin{equation}\label{eqModerate}
\frac{\codim W}{a_W}\le \frac{\codim W'}{a_{W'}}\quad\hbox{for any
two dense edges}\,\, W\subset W'.
\end{equation}

\begin{lem}\label{lemModerate} If $\cA$ is a hyperplane arrangement of moderate type, then
$(\codim W)/a_W$ is a jumping number of $\cA$ for any dense edge
$W$.
\end{lem}
\begin{proof} For any given dense edge $W$, we may assume that
there is strict inequality in (\ref{eqModerate}) by replacing $W$
with an edge containing it without changing $(\codim W)/a_W$. The
conclusion follows directly from the formula for multiplier ideals
in terms of dense edges of \cite{Te}.
\end{proof}

By \cite{ELSV}, Lemma \ref{lemModerate} implies:

\begin{thm}\label{corModerate} If $\cA$ is a hyperplane arrangement of moderate type, then
for any dense edge $W$, $-(\codim W)/a_W$ is a root of the
$b$-function of a polynomial defining $\cA$. So Conjecture
\ref{conjBernstein} holds in this case.
\end{thm}

\begin{rem}\label{remModerate} (i) Lemma \ref{lemModerate} does not hold if $\cA$ is not of moderate type
as is shown in Remark \ref{remCounterex}.

\medskip
(ii) It is also possible to prove Theorem \ref{corModerate} by
using \cite{Sa}, 4.2(b).

\medskip
(iii) If  $\cA$ is central, essential, and reduced, and $n=3$, then condition
(\ref{eqModerate}) holds if and only if the multiplicity along any
codimension 2 edge is at most $2d/3$, where $d$ is the number
of hyperplanes of $\cA$. Hence Theorem
\ref{corModerate} strengthens Example \ref{exdim3}.
\end{rem}

\end{document}